\documentstyle[11pt]{article}

\font\goth=eufm10

\font\bb=msbm10

\newtheorem{exemple}{Exemple}

\newtheorem{theoreme}[exemple]{Th\'eor\`eme}
\newtheorem{theorem}[exemple]{Theorem}
\newtheorem{corollaire}[exemple]{Corollaire}
\newtheorem{corollary}[exemple]{Corollary}

\newtheorem{conjecture}[exemple]{Conjecture}

\def\boxit#1#2{\setbox1=\hbox{\kern#1{#2}\kern#1}%
\dimen1=\ht1 \advance\dimen1 by #1 \dimen2=\dp1 \advance\dimen2 by #1
\setbox1=\hbox{\vrule height\dimen1 depth\dimen2\box1\vrule}%
\setbox1=\vbox{\hrule\box1\hrule}%
\advance\dimen1 by .4pt \ht1=\dimen1
\advance\dimen2 by .4pt \dp1=\dimen2 \box1\relax}

\def\adots{\mathinner{\mkern2mu\raise1pt\hbox{.}
\mkern3mu\raise4pt\hbox{.}\mkern1mu\raise7pt\hbox{.}}}
\def\<{\langle\,}
\def\>{\,\rangle}

\def\cf{{\it cf.$\ $}}

\def\End{{\rm End\, }}
\def\ev{{\rm ev}}

\def\N{\hbox{\bb N}}
\def\Z{\hbox{\bb Z}}
\def\C{\hbox{\bb C}}

\def\nn{\hbox{\goth n\hskip 1pt}}

\def\O{{\cal O}}
\def\slchap{\widehat{\hbox{\goth sl}}}
\def\H{\widehat H}
\def\m{{\bf m}}
\def\n{{\bf n}}
\def\p{{\bf p}}

\def\today{\number\day \space\ifcase\month\or
 Janvier\or F\'evrier\or Mars\or Avril\or Mai\or Juin\or
 Juillet\or Ao\^ut\or Septembre\or Octobre\or Novembre\or D\'ecembre\fi
 \space\number\year}

\title{{\bf Repr\'esentations induites d'alg\`ebres de Hecke \\ affines 
 et singularit\'es de $R$-matrices}} 

\author{{Bernard \sc  Leclerc} et {Jean-Yves \sc  Thibon}}

\date{}

\begin{document}

\newfont{\rml}{cmu10}
\newfont{\rmL}{cmr12}

\maketitle

\small
{\bf R\'esum\'e --- } Nous donnons un crit\`ere explicite 
d'irr\'eductibilit\'e pour des produits d'induction de modules
d'\'evaluation d'alg\`ebres de Hecke affines de type A.
Ce crit\`ere permet de d\'eterminer la forme des z\'eros
et des p\^oles de la $R$-matrice trigonom\'etrique associ\'ee
\`a une repr\'esentation d'\'evaluation  de
$U_v(\slchap_N)$.

\normalsize
\bigskip 

\centerline{{\bf Induced representations of affine Hecke algebras}}  
\centerline{{\bf and singularities of $R$-matrices}}

\bigskip
\small
{\it {\bf Abstract --- } We give an explicit criterion for the
irreducibility of some induction products of evaluation
modules of affine Hecke algebras of type $A$.
This allows to describe the form of the zeros and poles
of the trigonometric $R$-matrix associated to any evaluation 
module of $U_v(\slchap_N)$. }
\normalsize\rm


\vskip 12 mm
{\bf Abridged english version} --- \small
Let $\H_n(u)$ be the affine Hecke algebra associated to $GL(n)$
and $H_n(u)$ its finite-dimensional subalgebra
of type $A_{n-1}$. 
We assume that the parameter $u\in\C^*$ is not a root of unity.
The simple $H_n(u)$-modules $S_\lambda$ are parametrized by partitions
$\lambda$ of $n$.
Let $S_\lambda(z)$ be the $\H_n(u)$-module
obtained from $S_\lambda$ by evaluation at $z\in\C^*$.
Write $\lambda = (\lambda_1,\ldots ,\lambda_r)$
and $\lambda'=(l_1,\ldots ,l_k)$ (the conjugate of $\lambda$).
Let 
$$
{\cal E}_\lambda = \{e =\lambda_i+l_j-i-j+1 , \ 1\le i \le r,\ 1\le j \le k \}
$$
denote the set of hook-lengths of $\lambda$ and set 
$
{\cal Z}_\lambda = \{u^{\pm e},\ e\in {\cal E}_\lambda \}.
$
Given finite-dimensional $\H_{n_i}(u)$-modules $M_i,\ i=1,2$,
we write $M_1 \odot M_2$ for the induction product of $M_1$
and $M_2$ (this is then an $\H_{n_1+n_2}(u)$-module).
\begin{theorem}  Let $z_1,\ldots ,z_m \in \C^*$.
The module $S_\lambda(z_1)\odot \cdots \odot S_\lambda(z_m)$
is simple if and only~if 
${z_j/z_i} \not \in {\cal Z}_\lambda$ for all $i,\,j$.
\end{theorem}
The proof of Theorem~\ref{TH2} is reduced to a problem of canonical bases
in the following way.
Let $R_n$ be the complexified Grothendieck group of the category ${\cal C}_n$
of $\H_n(u)$-modules for which the generators $y_i$ of the maximal 
commutative subalgebra have eigenvalues of the form $u^k,\ k\in\Z$.
By Zelevinsky's classification \cite{Zel80}, 
the simple modules $L_\m$ of ${\cal C}_n$ are parametrized by 
the multisegments $\m = \sum_{i\le j} m_{ij} [i,j]$
over $\Z$ such that $\sum_{i<j} m_{ij}(j-i+1) = n$.
Let $R=\bigoplus_{n\ge 0} R_n$.
Following Zelevinsky we consider $R$ as a bialgebra
with multiplication and comultiplication corresponding 
respectively to induction and restriction with respect
to the maximal parabolic subalgebras $\H_k(u)\otimes\H_{n-k}(u)$.
Zelevinsky has shown that $R$ is a polynomial ring in
the generators $[L_{[i,j]}],\ i\le j$ associated with segments
$[i,j]$.
We denote by $M_\m$ the standard induced module corresponding to $\m$,
so that $[M_\m]=\prod_{i\le j}[L_{[i,j]}]^{m_{ij}}$.

On the other hand, let $A^-=A[N^-_\infty]$ denote the ring
of polynomial functions on the group
$N^-_\infty$ of lower unitriangular $\Z\times\Z$-matrices
with a finite number of non-zero off-diagonal entries.
It is a polynomial ring in the coordinate functions
$t_{ji},\ i<j$, and has a natural bialgebra structure.
We put $t_\m := \prod_{i\le j} t_{j+1,i}^{m_{ij}}$.
It is well known that $A^-$ and the enveloping algebra $U^-=U(\nn^-_\infty)$
are dual as bialgebras, the dual basis of $\{t_\m\}$ being
a basis $\{E_\m\}$ of Poincar\'e-Birkhoff-Witt type. 
Let $\{G(\m)\}$ be the canonical basis of $U^-$
obtained by specializing at $q=1$ the canonical basis of
$U_q(\nn^-_\infty)$ defined by Lusztig \cite{Lu90},
and let $\{G^*(\m)\}$ be the basis of $A^-$ dual to $\{G(\m)\}$.
It follows from the $p$-adic analogue of the Kazhdan-Lusztig
conjecture formulated by Zelevinsky \cite{Zel81} and proved
by Ginzburg \cite{CG}, and from the geometric description
by Lusztig \cite{Lu90} of the coefficients of the expansion of $G(\m)$
on $\{E_\m\}$ that
\begin{theorem}
The map $\Phi : R \longrightarrow A^-$ defined by 
$\Phi[M_\m]:=t_\m$ is an isomorphism of bialgebras,
and one has $\Phi[L_\m]=G^*(\m)$.
In particular, $L_\m\odot L_\n$ is simple if and only if $G^*(\m)G^*(\n)$
belongs to the dual canonical basis of $A^-$.
\end{theorem}
This theorem in a dual formulation is due to Ariki \cite{Ar}, 
but here we want to emphasize that when $u$ is not a root of unity, 
the dual canonical basis is more natural, as demonstrated below.

Let $A_q^-$ be the $q$-deformation of $A^-$, and let
$\{G_q^*(\m)\}$ denote the basis dual to the canonical basis
of $U_q(\nn^-_\infty)$.
The multiplicative properties of $\{G_q^*(\m)\}$ have
been studied by Berenstein and Zelevinsky.
They have formulated the following conjecture \cite{BZ}.
\begin{conjecture}
$$G_q^*(\m)\,G_q^*(\n)=q^k\,G_q^*(\p) \mbox{ for } k\in\Z
\ \ \Longleftrightarrow \ \
G_q^*(\m)\,G_q^*(\n)=q^lG_q^*(\n)\,G_q^*(\m) \mbox{ for } l\in\Z.
$$
\end{conjecture}
For $a\in\Z$, let $G^*(\lambda,a)=\Phi[S_\lambda(u^a)]$
be the element of the dual canonical basis corresponding
to the module $S_\lambda(z)$ evaluated at $z=u^a$.
The $G_q^*(\lambda,a)$ are distinguished elements of $\{G_q^*(\m)\}$, 
namely they are the quantum flag minors of the matrix 
${\bf T}_q=[t^{(q)}_{ji}]$ of the $q$-coordinate functions
$t^{(q)}_{ji}\in A^-_q$ \cite{BZ,LZ}.
Using the explicit
description given in \cite{LZ} of all pairs of $q$-commuting
quantum flag minors, we prove the following theorem, from which
Theorem~\ref{TH2} follows.
\begin{theorem}
{\rm (i)} \ Conjecture~\ref{CJBZ} is true when 
$G_q^*(\m)=G_q^*(\lambda,a)$ and $G_q^*(\n)=G_q^*(\mu,b)$
are quantum flag minors.

\noindent
{\rm (ii)} \ The product $G_q^*(\lambda,a_1)\cdots G_q^*(\lambda,a_m)$
belongs to the dual canonical basis (up to a power of $q$)
if and only if $|a_i-a_j| \not \in {\cal E}_\lambda$ for all $i,j$.
\end{theorem}
We note that, using the quantum affine Schur-Weyl duality between
$\H_n(u)$ and $U_v(\slchap_N)$ (with $u=v^2$) \cite{Che, GRV, ChPr96},
one can deduce from Theorem~\ref{TH2} the form of the singularities
of the trigonometric $R$-matrix $R_\lambda(z)$ associated to any 
simple evaluation module $V_\lambda(z)$ of $U_v(\slchap_N)$ \cite{Ji}.
\begin{corollary}
All singularities of $R_\lambda(z)$ are contained in the set ${\cal Z}_\lambda$.  
\end{corollary}

\setcounter{equation}{0}
\setcounter{exemple}{0}

\normalsize

\vskip 0.4cm
\noindent
---------------------------------
\vskip 0.6cm

Soit $\H_n(u)$ l'alg\`ebre de Hecke affine associ\'ee \`a $GL(n)$.
C'est l'alg\`ebre sur $\C$ engendr\'ee par des \'el\'ements
inversibles $y_1,\ldots ,y_n,T_1, \ldots , T_{n-1}$ soumis aux
relations
\begin{eqnarray*}
&&T_iT_{i+1}T_i=T_{i+1}T_iT_{i+1},\ \quad 1\le i\le m-2,\label{EQ_T1}\\
&&T_iT_j=T_jT_i,\qquad\qquad \qquad \vert i-j\vert>1,\label{EQ_T2}\\
&&(T_i-u)(T_i+1)=0,\qquad \ 1\le i\le m-1,\label{EQ_T3}\\
&&y_iy_j=y_jy_i, \qquad\qquad\qquad \ 1\le i,j\le m,\label{EQ_YY}\\
&&y_jT_i=T_iy_j,\qquad\qquad \quad \ \ \mbox{ $j\not= i,i+1$},\label{EQ_TY2}\\
&&T_iy_iT_i=u\,y_{i+1}, \qquad\qquad \ \ 1\le i\le m-1 .\label{EQ_Y}
\end{eqnarray*}
Nous supposons que $u\in\C^*$ n'est pas une racine de l'unit\'e.
Notons $H_n(u)$ la sous-alg\`ebre de $\H_n(u)$ engendr\'ee
par $T_1, \ldots , T_{n-1}$, et pour $z\in\C^*$ soit 
$\ev_z : \H_n(u) \longrightarrow H_n(u)$
l'homomorphisme d'\'evaluation d\'efini par $\ev_z(T_i) = T_i$
et $\ev_z(y_1) = z$. 
Les $H_n(u)$-modules simples $S_\lambda$ sont param\'etr\'es
par les partitions $\lambda$ de $n$.
Soit $S_\lambda(z)$ le $\H_n(u)$-module d\'eduit de $S_\lambda$
par \'evaluation en $z$.
Notons $\lambda = (\lambda_1,\ldots ,\lambda_r)$,
$\lambda'=(l_1,\ldots ,l_k)$ la partition conjugu\'ee de $\lambda$,
$$
{\cal E}_\lambda = \{e =\lambda_i+l_j-i-j+1 , \ 1\le i \le r,\ 1\le j \le k \}
$$
l'ensemble des longueurs d'\'equerres de $\lambda$, et posons 
$
{\cal Z}_\lambda = \{u^{\pm e},\ e\in {\cal E}_\lambda \}.
$
Etant donn\'es des $\H_{n_i}(u)$-modules
de dimension finie $M_i,\ i=1,2$, nous d\'esignons par 
\[
M_1 \odot M_2 := M_1 \otimes M_2 \uparrow_{\H_{n_1} \otimes \H_{n_2}}^{\H_{n_1+n_2}}
\]
le $\H_{n_1+n_2}(u)$-module obtenu par induction.
\begin{theoreme} \label{TH2} Soient $z_1,\ldots , z_m$ des complexes non nuls
et $\lambda$ une partition de $n$. 
Le produit d'induction $$S_\lambda(z_1)\odot \cdots \odot S_\lambda(z_m)$$
est un $\H_{nm}(u)$-module simple si et seulement si 
${z_j/z_i} \not \in {\cal Z}_\lambda$ pour tous $i,j$.
\end{theoreme}
Pour d\'emontrer le Th\'eor\`eme~\ref{TH2}, nous nous ramenons
\`a un probl\`eme de bases cano\-ni\-ques.
Soit $R_n$ le groupe de Grothendieck complexifi\'e de la cat\'egorie 
${\cal C}_n$ des
$\H_n(u)$-modules pour lesquels les g\'en\'erateurs $y_i$ 
de la sous-alg\`ebre commutative maximale ont toutes leurs
valeurs propres de la forme $u^k,\ k\in\Z$. 
D'apr\`es la classification de Zelevinsky \cite{Zel80}, 
les modules simples $L_\m$ de ${\cal C}_n$ sont param\'etr\'es
par les multi-segments $\m = \sum_{i\le j} m_{ij} [i,j]$
\`a support dans $\Z$ tels que $\sum_{i<j} m_{ij}(j-i+1) = n$.
On pose $R=\bigoplus_{n\ge 0} R_n$ (o\`u $R_0 = \C$).
Zelevinsky a muni $R$ d'une structure de big\`ebre,
la multiplication correspondant au produit d'induction
et la comultiplication aux diverses restrictions
de $\H_n(u)$ aux sous-alg\`ebres paraboliques $\H_k(u)\otimes\H_{n-k}(u)$.
Zelevinsky a montr\'e que $R$ est l'anneau des polyn\^omes en les
g\'en\'erateurs $[L_{[i,j]}],\ i\le j$ associ\'es aux segments
$[i,j]$, et que la comultiplication est donn\'ee sur ces g\'en\'erateurs
par
\begin{equation}\label{EQ1}
c[L_{[i,j]}] = 1\otimes [L_{[i,j]}] + 
\sum_{k=i}^{j-1}\  [L_{[i,k]}]\otimes [L_{[k+1,j]}]
+ [L_{[i,j]}]\otimes 1 \,.
\end{equation}
Zelevinsky a d\'efini un ordre partiel $\unlhd$
sur l'ensemble des multi-segments et a montr\'e que les
facteurs de composition du module induit standard $M_\m$
associ\'e au multisegment $\m$ \'etaient les $L_\n$
avec $\m\unlhd\n$. On a donc dans $R$
\begin{equation}\label{EQ2}
[M_\m] = \sum_{\hbox{\footnotesize$\m\unlhd\n$}} K_{\m\n} [L_\n]\,,
\end{equation}
pour certains entiers positifs $K_{\m\n}$.
Les multi-segments param\`etrent \'egalement les orbites
nilpotentes gradu\'ees $\O_\m$, (c'est-\`a-dire les 
classes d'isomorphismes de repr\'esentations du carquois
de type $A_\infty$),
et on a $\overline{\O}_\n = \coprod_{\hbox{\footnotesize$\m\unlhd\n$}} \O_\m$.
Zelevinsky \cite{Zel81} a formul\'e un analogue $p$-adique
de la conjecture de Kazhdan-Lusztig :
\begin{equation}\label{EQ3}
K_{\m\n}  = \sum_{i\ge 0} \dim {\cal H}^i(\overline{\O}_\n)_\m \,,
\end{equation}
o\`u ${\cal H}^i(\overline{O}_\n)_\m$ d\'esigne la tige en un
point $x\in \O_\m$ du $i$-\`eme faisceau de cohomologie 
d'intersection de la vari\'et\'e $\overline{\O}_\n$.
Cette conjecture a \'et\'e d\'emontr\'ee par Ginzburg
(\cf \cite{CG}, Theorem 8.6.23).

Par ailleurs, soit $A^-=A[N^-_\infty]$ l'anneau des fonctions
polynomiales sur le groupe $N^-_\infty$ des $\Z\times\Z$-matrices
unitriangulaires inf\'erieures ayant un nombre fini de coefficients
non-nuls hors de la diagonale. C'est l'anneau de polyn\^omes
en les fonctions coordonn\'ees $t_{ji},\ i<j$.
La comultiplication naturelle est donn\'ee sur ces g\'en\'erateurs
par 
\begin{equation}\label{EQ4}
\delta t_{ji} = t_{ji}\otimes 1 +\sum_{k=i+1}^{j-1} t_{jk}\otimes t_{ki}
+ 1\otimes t_{ji}\,.
\end{equation}
Au multi-segment $\m = \sum_{i\le j} m_{ij} [i,j]$ nous
associons le mon\^ome $t_\m := \prod_{i\le j} t_{j+1,i}^{m_{ij}}$.
On sait que la big\`ebre $A^-$ est duale de l'alg\`ebre 
enveloppante $U(\nn^-_\infty)$ de l'alg\`ebre de Lie de~$N^-_\infty$.
La base duale de $\{t_\m\}$ est la base de Poincar\'e-Birkhoff-Witt
$E_\m:= \overrightarrow{\prod_{[i,j]}}\, 
E_{j+1,i}^{m_{ij}}/m_{ij}!$,
o\`u les $E_{j+1,i}$ sont les unit\'es matricielles, et les
segments $[i,j]$ sont ordonn\'es par 
$$
[i,j] < [k,l] \quad \Longleftrightarrow \quad
(j<l 
\mbox{ ou }
\mbox{($j=l$ et $i<k$)).}
$$
Soit $\{G(\m)\}$ la base canonique de $U(\nn^-_\infty)$ obtenue
en sp\'ecialisant \`a $q=1$ la base canonique de $U_q(\nn^-_\infty)$
d\'efinie par Lusztig \cite{Lu90}. 
D'apr\`es l'\'equation (\ref{EQ3}) et \cite{Lu90}, 10.7, on~a 
\begin{equation}\label{EQ5}
G(\n) = \sum_{\hbox{\footnotesize$\m\unlhd\n$}} K_{\m\n} \, E_\m\,.
\end{equation}
Notons $\{G^*(\m)\}$ la base de $A^-$ duale de $\{G(\m)\}$.
Il r\'esulte de (\ref{EQ1}) (\ref{EQ2}) (\ref{EQ4}) (\ref{EQ5}):
\begin{theoreme}
L'application $\Phi : R \longrightarrow A^-$ d\'efinie par
$\Phi[M_\m]:=t_\m$ est un isomorphisme de big\`ebres, 
et on a $\Phi[L_\m]=G^*(\m)$.
En particulier, le produit d'induction $L_\m\odot L_\n$
est simple si et seulement si le produit $G^*(\m)G^*(\n)$
appartient \`a la base canonique duale de $A^-$.
\end{theoreme}
Ce th\'eor\`eme, sous une forme duale, est d\^u \`a Ariki \cite{Ar}, 
mais nous voudrions faire remarquer ici que lorsque $u$
n'est pas une racine de l'unit\'e, la base canonique duale
est plus naturelle, ainsi qu'on va le voir maintenant.

Soit $A_q^-$ la $q$-d\'eformation de $A^-$, et notons
$\{G_q^*(\m)\}$ la base duale de la base canonique de $U_q(\nn^-_\infty)$.
D'apr\`es \cite{Lu91} 11.5 (b), si l'on \'ecrit 
\begin{equation}
G_q^*(\m)\,G_q^*(\n)=\sum_\p \alpha_{\m\n}^\p(q)\,G_q^*(\p)\,,
\end{equation}
on a $\alpha_{\m\n}^\p(q) \in \N[q,q^{-1}]$.
Donc 
\[
G^*(\m)\,G^*(\n)=G^*(\p) \quad \Longleftrightarrow \quad
G_q^*(\m)\,G_q^*(\n)=q^k\,G_q^*(\p) \mbox{ pour } k\in\Z.
\]
Les propri\'et\'es multiplicatives de $\{G_q^*(\m)\}$ ont \'et\'e 
\'etudi\'ees par Berenstein et Zelevinsky.
Ils ont en particulier formul\'e la conjecture suivante \cite{BZ}.
\begin{conjecture}\label{CJBZ}
$$G_q^*(\m)\,G_q^*(\n)=q^k\,G_q^*(\p) \mbox{ pour } k\in\Z
\ \ \Longleftrightarrow \ \  
G_q^*(\m)\,G_q^*(\n)=q^lG_q^*(\n)\,G_q^*(\m) \mbox{ pour } l\in\Z.
$$
\end{conjecture}
Pour $a\in\Z$, notons $G^*(\lambda,a)=\Phi[S_\lambda(u^a)]$
les vecteurs de la base canonique duale cor\-res\-pon\-dant aux modules
d'\'evaluation en $z=u^a$.
Les $G_q^*(\lambda,a)$ forment un sous-ensemble distingu\'e
de la base canonique duale : ce sont les mineurs quantiques drapeaux
de la matrice ${\bf T}_q=[t^{(q)}_{ji}]$ des $q$-fonctions coordonn\'ees
$t^{(q)}_{ji}\in A^-_q$ \cite{BZ,LZ}.
En nous appuyant sur la description
explicite donn\'ee dans \cite{LZ} de toutes les paires de mineurs
quantiques drapeaux qui commutent \`a une puissance de $q$ pr\`es,
nous obtenons le th\'eor\`eme suivant qui entra\^\i ne le Th\'eor\`eme~\ref{TH2}.
\begin{theoreme}\label{TH5}
{\rm (i)} \ La Conjecture~\ref{CJBZ} est vraie lorsque 
$G_q^*(\m)=G_q^*(\lambda,a)$ et $G_q^*(\n)=G_q^*(\mu,b)$
sont des mineurs quantiques drapeaux.

\noindent
{\rm (ii)}  \ Le produit $G_q^*(\lambda,a_1)\cdots G_q^*(\lambda,a_m)$
appartient \`a la base canonique duale (\`a une puissance de $q$ pr\`es)
si et seulement si $|a_i-a_j| \not \in {\cal E}_\lambda$ pour tous $i,j$.
\end{theoreme} 
\medskip
\noindent
{\bf Remarque } Il est raisonnable de conjecturer que $L_\m \odot L_\n$
est irr\'eductible si et seulement si $L_\m\odot L_\n$ est
isomorphe \`a $L_\n\odot L_\m$ (\cf \cite{AK}, \cite{NT}).
Il serait int\'eressant d'\'eclaircir le lien entre cette
conjecture et la Conjecture~\ref{CJBZ} de Berenstein et Zelevinsky.

\medskip 
Soit $U_v(\slchap_N)$ l'alg\`ebre affine quantique de param\`etre 
$v=u^{1/2}$.
Soit $\lambda$ une partition de longueur $\le N$.
L'homomorphisme d'\'evaluation en $z\in\C^*$
de $U_v(\slchap_N)$ vers $U_v(gl_N)$
permet de munir le $U_v(gl_N)$-module simple $V_\lambda$ de plus haut poids 
$\lambda$  d'une structure de $U_v(\slchap_N)$-module de dimension
finie, not\'e $V_\lambda(z)$ \cite{Ji}.
En utilisant la dualit\'e de Schur-Weyl affine quantique \cite{Che, GRV, ChPr96},
on d\'eduit du Th\'eor\`eme~\ref{TH2} que si pour tous $i,j$
on a $z_j/z_i \not \in {\cal Z}_\lambda$, alors le produit tensoriel
$V_\lambda(z_1)\otimes \cdots \otimes V_\lambda(z_m)$ est un
$U_v(\slchap_N)$-module irr\'eductible.
Soit $R_\lambda(z)$ la $R$-matrice trigonom\'etrique associ\'ee
\`a $V_\lambda(z)$ \cite{Ji}.
 C'est une fonction rationnelle de $z$
\`a valeurs dans $\End(V_\lambda\otimes V_\lambda)$, qui est
solution de l'\'equation de Yang-Baxter quantique
\[
R_\lambda^{12}(z_1/z_2) \, R_\lambda^{13}(z_1/z_3) \, R_\lambda^{23}(z_2/z_3)
=
R_\lambda^{23}(z_2/z_3) \, R_\lambda^{13}(z_1/z_3) \, R_\lambda^{12}(z_1/z_2)
\]
pour des nombres complexes g\'en\'eriques $z_1,\,z_2,\,z_3$.
\begin{corollaire}
Les singularit\'es de $R_\lambda(z)$ sont contenues dans l'ensemble ${\cal Z}_\lambda$.
\end{corollaire}

\bigskip
\noindent
{\small
{\bf Remerciements.} Ce travail est n\'e de conversations avec 
M. Nazarov, lors d'une visite de B.~L. \`a l'Universit\'e de York en 1998.
Il a \'et\'e continu\'e pendant notre s\'ejour au R.I.M.S. de Kyoto 
dans le cadre du projet 1998 {\it ``Combinatorial methods in representation
theory''} organis\'e par M.~Kashiwara, K. Koike, S. Okada, I. Terada, H.-F. Yamada.
Nous remercions A. Lascoux et M. Nazarov pour de pr\'ecieuses
et stimulantes discussions, et les Universit\'es de York et de Kyoto
pour leur accueil chaleureux.

}

\bigskip \scriptsize
\noindent
\begin{tabular}{ll}
{\sc B.  Leclerc} : &
D\'epartement de Math\'ematiques,
Universit\'e de Caen, Campus II,\\
& Bld Mar\'echal Juin,
BP 5186, 14032 Caen cedex, France\\
&email : {\tt leclerc@math.unicaen.fr}\\[5mm]
{\sc J.-Y. Thibon} :&
Institut Gaspard Monge,
Universit\'e de Marne-la-Vall\'ee,\\
&Champs-sur-Marne,
77454 Marne-la-Vall\'ee cedex 2, France\\
&email : {\tt jyt@weyl.univ-mlv.fr}
\end{tabular}

\end{document}